# Non-recursive equivalent of the conjugate gradient method without the need to restart

Josip Dvornik[1]   Damir Lazarevic[1,3]   Antonia Jaguljnjak Lazarevic[2]   Marija Demsic[1]

[1]University of Zagreb, Faculty of Civil Engineering, Kaciceva 26, Zagreb 10 000, Croatia

[2]University of Zagreb, Faculty of Mining, Geology and Petroleum Engineering, Pierottijeva 6, Zagreb 10 000, Croatia

## Abstract

A simple alternative to the conjugate gradient (CG) method is presented; this method is developed as a special case of the more general iterated Ritz method (IRM) for solving a system of linear equations. This novel algorithm is not based on conjugacy, i.e. it is not necessary to maintain overall orthogonalities between various vectors from distant steps. This method is more stable than CG, and restarting techniques are not required. As in CG, only one matrix-vector multiplication is required per step with appropriate transformations. The algorithm is easily explained by energy considerations without appealing to the **A**-orthogonality in $n$-dimensional space. Finally, relaxation factor and preconditioning-like techniques can be adopted easily.

## 1 Introduction

Let
$$\mathbf{A}\mathbf{x} = \mathbf{b} \qquad (1)$$
be a real linear system with a symmetric positive definite (SPD) matrix of order $n$. By IRM, the solution is sought through successive minimisation of the corresponding energy functional, or the quadratic form
$$f(\mathbf{x}) = \frac{1}{2}\mathbf{x}^\mathrm{T}\mathbf{A}\mathbf{x} - \mathbf{x}^\mathrm{T}\mathbf{b} \qquad (2)$$
inside a small subspace formed at each iteration step [1]. After the convergence criterion is reached, a solution is found that is close to the unique minimiser of $f(\mathbf{x})$. Geometrically, this is the point close to the centre of the hyper-ellipsoids $f(\mathbf{x}) = c$, where $c$ are arbitrary real constants.

## 2 Briefly about IRM

The main idea here is to present the solution increment by the discretised Ritz method:
$$\mathbf{p}_{(i)} = \mathbf{\Phi}_{(i)}\mathbf{a}_{(i)} \qquad (3)$$
where $\mathbf{\Phi}_{(i)} = [\,\boldsymbol{\phi}_{1,(i)}\;\;\boldsymbol{\phi}_{2,(i)}\;\;\ldots\;\;\boldsymbol{\phi}_{m,(i)}\,]$ is a matrix of linearly independent coordinate vectors, and $\mathbf{a}_{(i)}$ is the vector of corresponding coefficients. The energy decrement associated with (3) can be also expressed as the quadratic function
$$\Delta f(\mathbf{a}_{(i)}) = \frac{1}{2}\mathbf{a}_{(i)}^\mathrm{T}\,\overline{\mathbf{A}}_{(i)}\mathbf{a}_{(i)} - \mathbf{a}_{(i)}^\mathrm{T}\bar{\mathbf{r}}_{(i)} \qquad (4)$$

---

[3] Corresponding author: e-mail: `damir@grad.hr`, 0000-0002-7439-719X, tel.: +385 1 4639 640

where $\overline{\mathbf{A}}_{(i)} = \mathbf{\Phi}_{(i)}^{\mathrm{T}} \mathbf{A} \mathbf{\Phi}_{(i)}$ and $\overline{\mathbf{r}}_{(i)} = \mathbf{\Phi}_{(i)}^{\mathrm{T}} \mathbf{r}_{(i)}$ are the SPD generalised (Ritz) matrix and the generalised residual vector, respectively, and both terms are of order $m$. After minimising (4), we obtain the system of equations that should be solved at each step:

$$\overline{\mathbf{A}}_{(i)} \mathbf{a}_{(i)} = \overline{\mathbf{r}}_{(i)} \qquad (5)$$

The solution is used to find the increment in (3), and $\mathbf{x}_{(i+1)} = \mathbf{x}_{(i)} + \mathbf{p}_{(i)}$ is updated afterwards. The residual is defined in a standard manner as $\mathbf{r}_{(i+1)} = \mathbf{b} - \mathbf{A}\mathbf{x}_{(i+1)}$.

Obviously, IRM represents an iterative procedure, where a discrete Ritz method is applied at each step, and a suitable set of coordinate vectors which span a subspace are generated. A local energy minimum is sought within that subspace (therefore (5) should be solved at each step), thereby decreasing the total energy of the system, which eventually converges to the required minimum. The subspace dimension, or the size of (5), is not limited. Rather, it aims to be small–much smaller than the number of unknowns ($m \ll n$), because every iteration must be as fast as possible. Such a small system (though $\overline{\mathbf{A}}_{(i)}$ is usually full) can be solved by any direct method. Simple pseudocode, with input data and sequence of instructions common for the iterative solution methods, is given by the Algorithm 1.

---

**Algorithm 1** Basic IRM algorithm

---

**Require:** $\mathbf{A}$, $\mathbf{b}$, $\mathbf{x}_{(0)}$, $\varepsilon$, $n_{max}$ 	{usually $\mathbf{x}_{(0)} \leftarrow \mathbf{0}$}

**Ensure:** $\mathbf{x}_{(i+1)}$ {close to $\mathbf{x}$}

1: $i \leftarrow 0$ 	{initialisation: steepest descent}

2: $\mathbf{r}_{(0)} \leftarrow \mathbf{b} - \mathbf{A}\mathbf{x}_{(0)}$

3: $q \leftarrow \mathbf{r}_{(0)}^{\mathrm{T}} \mathbf{r}_{(0)} / (\mathbf{r}_{(0)}^{\mathrm{T}} \mathbf{A} \mathbf{r}_{(0)})$

4: $\mathbf{p}_{(0)} \leftarrow q \mathbf{r}_{(0)}$

5: **while** $(\|\mathbf{r}_{(i)}\|_2 > \varepsilon \|\mathbf{r}_{(0)}\|_2) \wedge (i \leq n_{\max})$ **do** {Iterated Ritz method}

6:     $\mathbf{x}_{(i+1)} \leftarrow \mathbf{x}_{(i)} + \mathbf{p}_{(i)}$

7:     $\mathbf{r}_{(i+1)} \leftarrow \mathbf{b} - \mathbf{A}\mathbf{x}_{(i+1)}$

8:     generate $[\,\boldsymbol{\phi}_{1,(i)}\ \ \boldsymbol{\phi}_{2,(i)}\ \ \cdots\ \ \boldsymbol{\phi}_{m,(i)}\,]$

9:     $\overline{\mathbf{A}}_{(i)} \leftarrow [\,\boldsymbol{\phi}_{1,(i)}\ \ \boldsymbol{\phi}_{2,(i)}\ \ \cdots\ \ \boldsymbol{\phi}_{m,(i)}\,]^{\mathrm{T}} \mathbf{A} [\,\boldsymbol{\phi}_{1,(i)}\ \ \boldsymbol{\phi}_{2,(i)}\ \ \cdots\ \ \boldsymbol{\phi}_{m,(i)}\,]$

10:     $\overline{\mathbf{r}}_{(i)} \leftarrow [\,\boldsymbol{\phi}_{1,(i)}\ \ \boldsymbol{\phi}_{2,(i)}\ \ \cdots\ \ \boldsymbol{\phi}_{m,(i)}\,]^{\mathrm{T}} \mathbf{r}_{(i+1)}$

11:     $\mathbf{a}_{(i)} \leftarrow \overline{\mathbf{A}}_{(i)}^{-1} \overline{\mathbf{r}}_{(i)}$

12:     $\mathbf{p}_{(i+1)} \leftarrow [\,\boldsymbol{\phi}_{1,(i)}\ \ \boldsymbol{\phi}_{2,(i)}\ \ \cdots\ \ \boldsymbol{\phi}_{m,(i)}\,] \mathbf{a}_{(i)}$

13:     $i \leftarrow i + 1$

14: **end while** 	{end Iterated Ritz method}

---

The central problem involves quickly generating a small and efficient subspace, such that the energy reduction per step is as large as possible and the number of steps is extremely reduced. Usually, one coordinate vector is $\mathbf{p}_{(i)}$, and others are generated as $\mathbf{P}_j \mathbf{r}_{(i+1)}$, where $\mathbf{P}_j$ is fast approximation of $\mathbf{A}^{-1}$. Non-residual based generation ideas are also possible [2]. Because many strategies to construct $\mathbf{P}_j$ (or generally $\boldsymbol{\phi}_{j,(i)}$) exist, the algorithm also allows for any new routine that generate subspaces (e.g. those suggested by other independent researchers) to be easily implemented in (line 8, Alg. 1). Potentially, this may make the method even faster.

It should be noted that the conjugacy property is not explicitly taken into account in IRM, and coordinate vectors may become (almost) linearly dependent. Therefore, routines for subspace generation which prevent such a scenario are preferred and may even change between steps. Nevertheless, if this dependence arises, some pivots approach zero during the decomposition of $\overline{\mathbf{A}}_{(i)}$, which can be recognised and used for discarding corresponding equations from the small system. The subspace dimension is reduced in such cases, but $\overline{\mathbf{A}}_{(i)}$ becomes more regular and better conditioned. This strategy is faster than orthogonalisation, rejection, or replacement of dependent vectors [3].

IRM can also be considered as a generalisation of some iterative methods [1,2]. Depending on the choice of coordinate vectors, some solvers can be represented or interpreted as special cases of this approach. Furthermore, it is possible to combine good properties of several methods simultaneously. If appropriate vectors are selected, convergence should proceed faster than using any single method considered. Here, an improved CG algorithm (IRM-CG) is presented due to the popularity of SPD systems.

## 3  IRM-CG as a special case of IRM

The algorithm presented here also starts with the steepest descent (SD) step. Other steps are executed using a CG-like algorithm simulated by IRM with two coordinate vectors. The first vector is the current residual $\mathbf{r}_{(i+1)}$, and the second vector is previous solution increment $\mathbf{p}_{(i)}$. Vectors span a two-dimensional subspace. At each step, a system of two equations is solved and a new energy minimum within that plane is found (Algorithm 2).

This approach has three matrix-vector multiplications per step: one in line 7 and two in line 8. Applying two 'induced' recursive relations ('inherent' $\mathbf{A}$-orthogonalisation is not exploited here), only one such multiplication remains (as in CG). If line 11 (Alg. 2) is multiplied by $\mathbf{A}$, then

$$\mathbf{A}\mathbf{p}_{(i+1)} = [\,\mathbf{A}\mathbf{r}_{(i+1)} \quad \mathbf{A}\mathbf{p}_{(i)}\,]\,\mathbf{a}_{(i)} \tag{6}$$

Substituting $\boldsymbol{\alpha}_{(i)} = \mathbf{A}\mathbf{r}_{(i+1)}$ and $\boldsymbol{\beta}_{(i)} = \mathbf{A}\mathbf{p}_{(i)}$ into the first recursion yields

$$\boldsymbol{\beta}_{(i+1)} = [\,\boldsymbol{\alpha}_{(i)} \quad \boldsymbol{\beta}_{(i)}\,]\mathbf{a}_{(i)} \tag{7}$$

Second, the frequently used residual recursion $\mathbf{r}_{(i+1)} = \mathbf{r}_{(i)} - \mathbf{A}\mathbf{p}_{(i)}$ becomes

$$\mathbf{r}_{(i+1)} = \mathbf{r}_{(i)} - \boldsymbol{\beta}_{(i)} \tag{8}$$

Now, after the line 4 (Alg. 2), the new initialisation

$$\boldsymbol{\beta}_{(0)} \leftarrow \mathbf{A}\mathbf{p}_{(0)} \tag{9}$$

should be inserted, and the pseudocode inside the **while** loop becomes:

---

**Algorithm 2** Basic IRM-CG algorithm

---

**Require:** $\mathbf{A}$, $\mathbf{b}$, $\mathbf{x}_{(0)}$, $\varepsilon$, $n_{max}$ {usually $\mathbf{x}_{(0)} \leftarrow \mathbf{0}$}

**Ensure:** $\mathbf{x}_{(i+1)}$ {close to $\mathbf{x}$}

1: $i \leftarrow 0$ {initialisation: steepest descent}

2: $\mathbf{r}_{(0)} \leftarrow \mathbf{b} - \mathbf{A}\mathbf{x}_{(0)}$

3: $q \leftarrow \mathbf{r}_{(0)}^T \mathbf{r}_{(0)} / (\mathbf{r}_{(0)}^T \mathbf{A} \mathbf{r}_{(0)})$

4: $\mathbf{p}_{(0)} \leftarrow q\mathbf{r}_{(0)}$

5: **while** $(\|\mathbf{r}_{(i)}\|_2 > \varepsilon \|\mathbf{r}_{(0)}\|_2) \wedge (i \leq n_{\max})$ **do** {IRM-CG method}

6: $\quad \mathbf{x}_{(i+1)} \leftarrow \mathbf{x}_{(i)} + \mathbf{p}_{(i)}$

7: $\quad \mathbf{r}_{(i+1)} \leftarrow \mathbf{b} - \mathbf{A}\mathbf{x}_{(i+1)}$

8: $\quad \overline{\mathbf{A}}_{(i)} \leftarrow [\mathbf{r}_{(i+1)} \ \mathbf{p}_{(i)}]^T \mathbf{A} [\mathbf{r}_{(i+1)} \ \mathbf{p}_{(i)}]$

9: $\quad \overline{\mathbf{r}}_{(i)} \leftarrow [\mathbf{r}_{(i+1)}^T \mathbf{r}_{(i+1)} \ 0]^T$ {second term is zero because $\mathbf{r}_{(i+1)}^T \mathbf{p}_{(i)} = 0$}

10: $\quad \mathbf{a}_{(i)} \leftarrow \overline{\mathbf{A}}_{(i)}^{-1} \overline{\mathbf{r}}_{(i)}$

11: $\quad \mathbf{p}_{(i+1)} \leftarrow [\mathbf{r}_{(i+1)} \ \mathbf{p}_{(i)}] \mathbf{a}_{(i)}$

12: $\quad i \leftarrow i + 1$

13: **end while** {end IRM-CG method}

---

$\mathbf{x}_{(i+1)} \leftarrow \mathbf{x}_{(i)} + \mathbf{p}_{(i)}$

$\mathbf{r}_{(i+1)} \leftarrow \mathbf{r}_{(i)} - \boldsymbol{\beta}_{(i)}$

$\boldsymbol{\alpha}_{(i)} \leftarrow \mathbf{A}\mathbf{r}_{(i+1)}$ {sole matrix-vector multiplication}

$\overline{\mathbf{A}}_{(i)} \leftarrow [\mathbf{r}_{(i+1)} \ \mathbf{p}_{(i)}]^T [\boldsymbol{\alpha}_{(i)} \ \boldsymbol{\beta}_{(i)}]$ {$\overline{\mathbf{A}}_{(i)}$ is symmetric: $\mathbf{r}_{(i+1)}^T \boldsymbol{\beta}_{(i)} = \mathbf{p}_{(i)}^T \boldsymbol{\alpha}_{(i)}$} (10)

lines 9–11 remain unchanged

$\boldsymbol{\beta}_{(i+1)} \leftarrow [\boldsymbol{\alpha}_{(i)} \ \boldsymbol{\beta}_{(i)}] \mathbf{a}_{(i)}$

$i \leftarrow i + 1$

Due to round-off errors, as in CG, the residual is periodically (after $i_{\max}$ steps) updated from the equilibrium equation, i.e. the following should be used instead of the second line from (10):

**if** $i \bmod i_{\max} \neq 0$ **then**

  $\mathbf{r}_{(i+1)} \leftarrow \mathbf{r}_{(i)} - \boldsymbol{\beta}_{(i)}$

**else** (11)

  $\mathbf{r}_{(i+1)} \leftarrow \mathbf{b} - \mathbf{A}\mathbf{x}_{(i+1)}$

**endif**

## 4 Equivalence between CG and IRM-CG

The proof of equivalence between CG and IRM-CG is very simple, so it will be discussed only briefly. Initialisation is practically identical for both methods. In other steps, the minimum of the energy function inside the plane spanned by $\mathbf{r}_{(i+1)}$ and $\mathbf{p}_{(i)}$ is determined. In CG this is realised by $\mathbf{A}$-orthogonalisation and by solving a linear system of two equations in IRM-CG.

If exact arithmetic is considered, IRM-CG and CG have an identical sequence of intermediate results. The exact solution is obtained after $m$ steps, where $m$ is the number of different 'active' eigenvalues. If $\mathbf{b}$ is represented as a sum of eigenvectors $\boldsymbol{\varphi}_j$, i.e. $\mathbf{b} = \sum a_j \boldsymbol{\varphi}_j$, eigenvectors (and their corresponding eigenvalues) with $a_j \neq 0$ may be called 'active' (or 'inactive' otherwise). Of course, $m$ can be found only if all $n$ eigenpairs are detected. Multiple eigenvalues should be counted as one, and 'inactive' eigenvalues are not counted at all. This comment is only important for theoretical considerations, as IRM-CG is interesting as an iterative, not direct solution method.

## 5 Advantages of IRM-CG over CG

During real calculations (with round-off errors), IRM-CG is more stable and behaves better than CG. First, restarting of this algorithm is not needed because $\mathbf{A}$-orthogonality is not exploited. Namely, error in $\mathbf{A}$-orthogonality also exists if the IRM-CG formulation is used, but it is not accumulated during calculation. Therefore, inherited errors from the $\mathbf{A}$-orthogonality decrease, although non-exact arithmetic (as in every numerical process) causes new errors and affects convergence.

Consider simple example with diagonal $\mathbf{A}$ ($a_{1,1} = 1$ and $a_{2,2} = \kappa$, where $\kappa$ is a condition number of $\mathbf{A}$). If $\mathbf{b} = \begin{bmatrix} 1 & 1 \end{bmatrix}^{\mathrm{T}}$, using rational arithmetic exact solution $\mathbf{x} = \begin{bmatrix} 1 & 1/\kappa \end{bmatrix}^{\mathrm{T}}$ is obtained in two steps by both methods. To check stability of the methods, after the initialisation phase (for $i = 0$) a small disturbance $\delta$ to the second term of $\mathbf{p}_{(1)}$ is added (Fig. 1a). The second step of IRM-CG still gives exact result, but CG gives only s perturbed approximate solution as a function of $\delta$ and $\kappa$:

$$\begin{aligned}
\widetilde{x}_1 &= 2\left(\frac{1}{1+\kappa} + \frac{2(\kappa-1)}{4(\kappa+1) - 4\delta(\kappa-1) + \delta^2(\kappa-1)^2}\right) \\
\widetilde{x}_2 &= 2\left(\frac{1}{1+\kappa} + \frac{(\kappa-1)(\delta(\kappa-1)-2)}{\kappa(4(\kappa+1) - 4\delta(\kappa-1) + \delta^2(\kappa-1)^2)}\right)
\end{aligned} \quad (12)$$

Notice complexity of the CG solution, even with a diagonal matrix of order two. For better explanation of the expressions, over domain $\delta \in [-10^{-2}, 10^{-2}]$ functions $\widetilde{\mathbf{x}} - \mathbf{x} = f(\delta)$ for three values of $\kappa$, and $\widetilde{\mathbf{x}} - \mathbf{x} = f(\delta, \kappa)$ for $\kappa = 10^a$ ($a \in [0, 4]$) are given (Figs. 1b and c). Only if $\delta = 0$ exact solution is recovered from (12). Also, functions are non-symmetric, so CG behaves differently for $\pm\delta$.

When large equation systems are considered, $\delta$ is accumulated primarily due to the loss of $\mathbf{A}$-orthogonality, which is inherited recursively. The main reason is approximate CG solution at each step (in the current plane); in practice more complicated then (12). On the contrary, IRM-CG finds numerically 'exact' solution at each plane. Therefore system of two equations is repeatedly solved. Roughly, if $\delta$ is split into two components at each step, the one inside and the other orthogonal to the plane, the first component is 'exactly' resolved and does not produce inherited error. In CG, both components cause propagation of error. In other words, if $\delta$ lies on the plane, the behaviour of IRM-CG is as if $\delta = 0$, but if it is orthogonal to the plane the solution is disturbed.

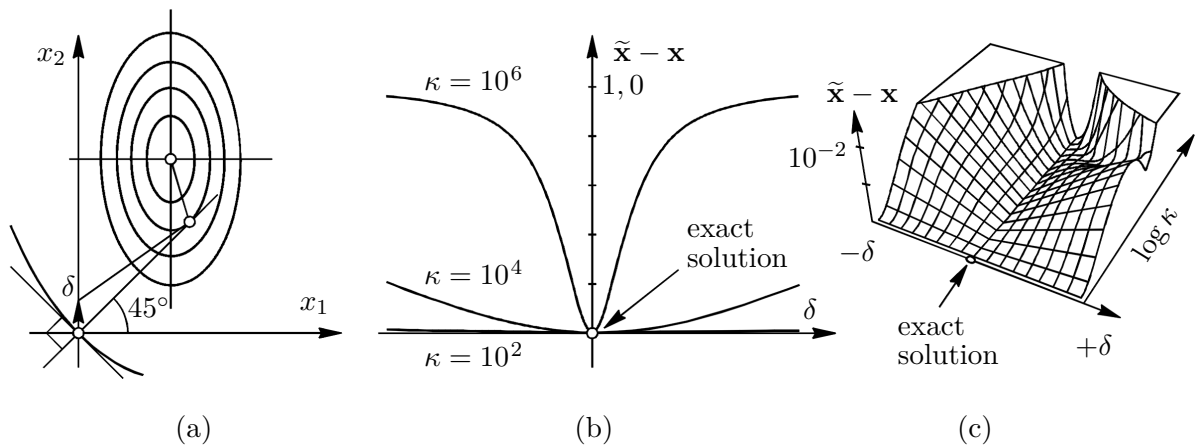

Figure 1: Stability of CG method: (a) interpretation of disturbance $\delta$, (b) function $\widetilde{\mathbf{x}} - \mathbf{x} = f(\delta)$, (c) surface $\widetilde{\mathbf{x}} - \mathbf{x} = f(\delta, \kappa)$
.

It is possible to interchange methods because the two approaches are equivalent. Each step may be executed by CG or IRM-CG, no matter how the earlier steps were performed. If CG is used as a solution method, it is suggested that one equivalent IRM-CG step be executed after some number of steps, but before orthogonality error becomes too large. This may be called "refresh" instead of traditionally "restart".

The second advantage of this formulation is the natural adoption of the relaxation factor $\omega \in (0, 2)$, known from the Successive overrelaxation method [4], and which may even change at each step. Line 6 (Alg. 2) is simply

$$\mathbf{x}_{(i+1)} \leftarrow \mathbf{x}_{(i)} + \omega_{(i)} \mathbf{p}_{(i)} \tag{13}$$

and the second term of $\overline{\mathbf{r}}_{(i)}$ is not zero (line 9, Alg. 2); it is rather $\omega_{(i)} \mathbf{r}_{(i+1)}^\mathrm{T} \mathbf{p}_{(i)}$. Also, recurrence relation in (11) becomes $\mathbf{r}_{(i+1)} \leftarrow \mathbf{r}_{(i)} - \omega_{(i)} \boldsymbol{\beta}_{(i)}$. Obviously, using $\omega \neq 1$, $\mathbf{A}$-orthogonality is lost, but convergence is improved in many cases [5]. The third advantage is the very natural adoption of (multi) preconditioning-like techniques [6]. In such cases, only line 8 (Alg. 2) is

reformulated as

$$\overline{\mathbf{A}} \leftarrow [\,\mathbf{p}_{(i)} \;\; \phi_{2,(i)} \;\; \phi_{3,(i)} \;\; \ldots\,]^{\mathrm{T}} \mathbf{A} \,[\,\mathbf{p}_{(i)} \;\; \phi_{2,(i)} \;\; \phi_{3,(i)} \;\; \ldots\,] \qquad (14)$$

The coordinate vectors are

$$\phi_{j,(i)} = \mathbf{M}_j^{-1} \mathbf{r}_{(i+1)}, \qquad 2 \leq j \leq n \qquad (15)$$

where $\mathbf{M}_j$ is a matrix according to the standard approach [7], which is used to produce a better conditioned system $\mathbf{M}_j^{-1}\mathbf{A}\mathbf{x} = \mathbf{M}_j^{-1}\mathbf{b}$ equivalent to (1). However, transformations of the CG algorithm required for such strategies are not needed here. According to IRM, that is just another way of generating coordinate vector(s). During the solution process, they can also become (exactly or approximately) linearly dependent, and one or several of them should be excluded.

Many possibilities to rapidly construct $\mathbf{M}^{-1}$ exist, such as (not always robust) incomplete Cholesky factorisations with different fill-ins [8], algebraic multigrid methods [9,10], and sparse approximate inverses [11]. It is even possible to use methods that are not useful as standalone solvers, as they are neither convergent nor numerically stable. For smoothing purposes, forward and backward techniques or any other promising order of unknowns may be useful.

## 6 Illustrative example

Consider a simple linear FEM benchmark: a cube discretised by 8-node solid elements, supported by the corner springs of stiffnesses $k$, and loaded with a vertical force at the top. This model has 3 993 unknowns (Fig. 2). The condition number is $\kappa \approx 8 \cdot 10^4$, which is calculated as the ratio of extreme eigenvalues. CG and IRM-CG behave almost identically (curves practically collide) for such a well-conditioned system. If the spring stiffness A is greatly reduced to $10^{-10}k$, then $\kappa \approx 3 \cdot 10^{13}$ and IRM-CG behaves much better than CG. Processes are terminated after $\varepsilon = 10^{-10}$ is reached. Of course, CG may be improved by preconditioning and restarting techniques. However, IRM-CG may also be enhanced by using additional coordinate vectors, while restarting strategies are not needed at all, as previously mentioned.

## 7 Conclusion

Although general theorems and proofs about algorithm convergence rate and stability are not given here, according to the results of numerical experiments with exact and floating-point arithmetic, IRM–CG should be an interesting replacement for a standard or preconditioned CG. Recursive $\mathbf{A}$–orthogonalisation, restarting recommendations, and transformations (necessary for preconditioning) are not required, hence the method should be very useful for solving non-well-posed problems. Finally, the property of conjugacy, which underlies many iterative procedures and is valid only for linear systems, is not absolutely necessary here. Therefore, IRM–CG can be successfully applied to nonlinear problems (including optimisation), where conjugacy is not even defined. This issue is vitally important, as iterative methods are used exclusively in these cases [12].

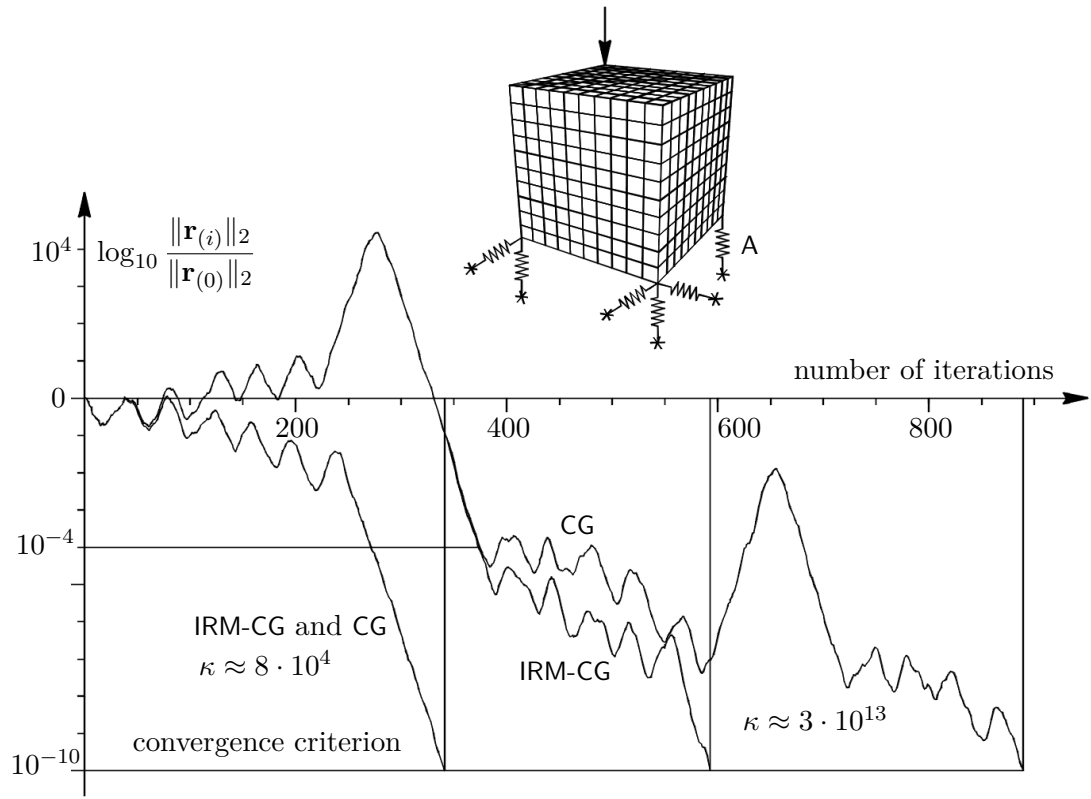

Figure 2: Behaviour of CG and IRM-CG for well-posed and ill-posed problems

## Acknowledgement

This work was fully supported by the Croatian Science Foundation under the project IP–2014–09–2899.

Non-recursive equivalent of the conjugate gradient method without the need to restart  9